\documentclass[12pt, a4paper]{article}
\title{On 3-matrix factorizations of polynomials}
\author{ Yves Baudelaire Fomatati$^{a}$\footnote{a. Department of Mathematics and Statistics, University of Ottawa, yfomatat@uottawa.ca .}}

\date{}

\usepackage[square,comma,numbers,sort&compress]{natbib}
\usepackage[centertags]{amsmath}
\usepackage{latexsym}
\usepackage{leqno}
\usepackage{amsfonts}
\usepackage{amsmath,amssymb}
\usepackage{enumerate}

\usepackage{amsthm}
\usepackage[utf8,latin9]{inputenc}

\theoremstyle{plain}
\newtheorem{remark}{Remark}[section]
\theoremstyle{plain}
\newtheorem{lemma}{Lemma}[section]
\theoremstyle{plain}
\newtheorem{proposition}{Proposition}[section]
\theoremstyle{plain}
\newtheorem{theorem}{Theorem}[section]
\theoremstyle{plain}
\newtheorem{definition}{Definition}[section]
\theoremstyle{plain}

\theoremstyle{plain}

\theoremstyle{plain}

\newtheorem{example}{Example}[section]

\theoremstyle{plain}

\usepackage[all, 2cell]{xy}
\usepackage{txfonts}

\begin{document}

\maketitle

\begin{quote}
  \textbf{Abstract}
\end{quote}
Let $R=K[x_{1},x_{2},\cdots, x_{m}]$ and $S=$ $K[y_{1},y_{2},\cdots, y_{m}]$ where $K$ is a field. 
In this paper, we propose a method showing how to obtain $3$-matrix factors for a given polynomial using either the Doolittle or the Crout decomposition techniques that we apply to matrices whose entries are not real numbers but polynomials.
We also define the category of $3$-matrix factorizations of a polynomial $f$ whose objects are $3$-matrix factorizations of $f$, that is triplets $(P,Q,T)$ of $m\times m $ matrices such that $PQT=fI_{m}$. Moreover, we construct a bifunctorial operation $\overline{\otimes}_{3}$ which is such that if $X$ (respectively $Y$) is a $3-$matrix factorization of $f\in R$ (respectively $g\in S$), then $X\overline{\otimes}_{3} Y$ is a $3-$matrix factorization of $fg\in K[x_{1},x_{2},\cdots, x_{m},y_{1},y_{2},\cdots, y_{m}]$. We call $\overline{\otimes}_{3}$ the multiplicative tensor product of $3-$matrix factorizations. Finally, we give some properties of the operation $\overline{\otimes}_{3}$.
In the sequel, $K$ is a field and except otherwise stated, our polynomials will be taken from $R=K[x_{1},x_{2},\cdots, x_{m}]$, where $K=\mathbb{R}$, the set of real numbers. Sometimes instead of indexing the indeterminates when they are at most three, we will write $x,y,z$.

\section{Introduction}

  Eisenbud was the first to introduce the notion of matrix factorization in his seminal paper \cite{eisenbud1980homological} in 1980. He proved that
matrix factorizations of $f\in K[[x]]$ describe all maximal Cohen-Macaulay modules (MCM modules)
without free summands (See \cite{leuschke2012cohen} for notes on MCM modules). This notion generalizes the classical polynomial factorization in the sense that classical polynomial factors can now be seen as $1\times 1$ matrix factors. The polynomial
$g=x^{3}+y^{2}$ is irreducible over $\mathbb{R}[x,y]$ but can be factorized as follows:
$$\begin{bmatrix}
    x  &  -y      \\
    y  &  x^{2}
\end{bmatrix}
\begin{bmatrix}
    x^{2}  &  y     \\
    -y  &   x
\end{bmatrix}
= (x^{3} + y^{2})\begin{bmatrix}
    1  &  0      \\
    0  &  1
\end{bmatrix}
=gI_{2} $$
$
(\begin{bmatrix}
    x  &  -y      \\
    y  &  x^{2}
\end{bmatrix},
\begin{bmatrix}
    x^{2}  &  y      \\
    -y  &  x
\end{bmatrix})
$
 is said to be a $2 \times 2$ matrix factorization of $g$.
 Matrix factorizations and some of their properties were studied in several papers including \cite{eisenbud1980homological}, \cite{carqueville2016adjunctions}, \cite{crisler2016matrix}, \cite{fomatati2022tensor} and \cite{camacho2015matrix}.
 In these papers, a matrix factorization of a polynomial $f$ is a pair of $m\times m$ matrices $(P,Q)$ such that $fI_{m}=PQ$. In this paper, we will refer to this type of matrix factorization as $2$-matrix factorization because we have two matrix factors. The category of $2$-matrix factorizations was defined in \cite{yoshino1998tensor}. Its objects are $2$-matrix factorizations. We will define the notion of $3$-matrix factorization of a polynomial $f$ (which has to do with fatorizing polynomials using three matrices) and we will show that a triplet $(P,Q,T)$ of $m\times m$ matrices whose product equals $fI_{m}$, is an object in the category of $3$-matrix factorizations of $f$ which will be constructed in the paper (cf. section \ref{subsec category of 3-matrix facto of f}). \\
One obvious reason for studying matrix factorizations and their properties is that originally, Eisenbud \cite{eisenbud1980homological} proved that
matrix factorizations of $f\in K[[x]]$ describe all maximal Cohen-Macaulay modules (MCM modules)
without free summands.
Another interesting reason for studying matrix factorizations is that
irreducible polynomials can be factorized using matrices.
Furthermore,
 Buchweitz et al. \cite{buchweitz1987cohen} found that matrix factorizations of polynomials (over the reals) of the form $f_{n}=x_{1}^{2}+\cdots + x_{n}^{2}$, for $n=1,2,4$ and $8$ are related to the existence of composition algebras over $\mathbb{R}$ of dimension $1,2,4$ and $8$ namely the complex numbers, the quaternions and the octonians.
More on the importance of matrix factorizations with references can be found in the introduction of \cite{fomatati2022tensor}.\\
Eisenbud (p.15 of \cite{eisenbud1980homological}) in 1980 originally defined a matrix factorization of an element $f$ in a ring $R$ (with unity) to be an ordered pair of maps of free $R-$modules $\phi: F\rightarrow G$ and $\psi: G \rightarrow F$ s.t., $\phi\psi=f\cdot 1_{G}$ and $\psi\phi=f\cdot 1_{F}$. \\
In 1998, Yoshino \cite{yoshino1998tensor} defined a matrix factorization of $f$ to be a pair of matrices $(P,Q)$ such that $fI=PQ$. Diveris and Crisler in 2016 used this definition (cf. Definition 1 of \cite{crisler2016matrix}) of Yoshino.  In this paper, we follow suit and we refer to this type of matrix factorization of a polynomial $f$ as a 2-matrix factorization of $f$ as already mentioned above. Next, we extend this definition to $n$-matrix factorizations of a polynomial, for $n=3$ and $n=4$. In their paper published in 2016, Carqueville and Murfet defined a matrix factorization using linear factorizations and $\mathbb{Z}_{2}$-graded modules (cf. p.8 of \cite{carqueville2016adjunctions}).
Detailed explanations on how the different definitions mentioned above are interrelated can be found in \cite{fomatati2023note}.\\
Properties of 2-matrix factorizations were used in \cite{crisler2016matrix} to give the minimal 2-matrix factorization for a polynomial which is the sum of squares of 8 monomials.
They were also used in chapter 6 of \cite{fomatati2019multiplicative} to give necessary conditions for the existence of a Morita Context in the bicategory of Landau-Ginzburg models. Moreover, one of the properties of 2-matrix factorizations was used to conclude that a polynomial admits more than one pair of 2-matrix factors (cf. proposition 2.2 of \cite{fomatati2019multiplicative}).\\
It is good to recall that there is an algorithm referred to as the standard method which produces 2-matrix factors of any given polynomial. Details about this algorithm can be found in \cite{crisler2016matrix}, \cite{fomatati2023refined}.
 We will use the Doolittle or the Crout matrix decomposition methods on matrices whose entries are polynomial and not just real numbers to obtain an $n$-matrix factorization of a given polynomial $f$ for $n=3$ and $n=4$ from any of its 2-matrix factorizations.\\
 Moreover, we propose a bifunctorial operation $\overline{\otimes}_{3}$ which is such that if $X$ (respectively $Y$) is a $3-$matrix factorization of $f\in R$ (respectively $g\in S$), then $X\overline{\otimes}_{3} Y$ is a $3-$matrix factorization of $fg\in K[x_{1},x_{2},\cdots, x_{m},y_{1},y_{2},\cdots, y_{m}]$. We call $\overline{\otimes}_{3}$ the multiplicative tensor product of $3-$matrix factorizations. The operation $\overline{\otimes}_{3}$ is the counterpart of the multiplicative tensor product of $2$-matrix factorizations that was defined in \cite{fomatati2023refined} and used in conjunction with the Yoshino tensor product of matrix factorizations to reduce the size of matrix factors on the class of summand reducible polynomials.
  Finally, we give some properties of the operation $\overline{\otimes}_{3}$.
\\
This paper is organized as follows: In the next section, we give some preliminaries. In section 3, after recalling the definition of $2$-matrix factorizations, we define what a $3$-matrix factorization of a polynomial is and give some examples. The category of $3$-matrix factorizations is defined in section 4. In section 5, we define the multiplicative tensor product of $3$-matrix factorizations denoted $\overline{\otimes}_{3}$ and we prove that it is a bifunctorial operation after giving some examples. Moreover, we give some of its properties. Finally, further problems are discussed in section 6.

\section{Preliminaries}
In this section, we give some preliminaries.
\subsection{Doolittle and Crout matrix decomposition methods} \label{Doolittle and Crout methods}

Here, we recall two matrix decomposition methods that are well known in the literature, namely Doolittle and Crout matrix decomposition methods.\\
First, we recall the following definition.
\begin{definition} $LU$-Factorization, Doolittle and Crout decompositions  \label{LU factorization}\\
 A nonsingular matrix $M$ has an $LU$-factorization if it can be expressed as the product of a lower-triangular matrix $L$ and an upper triangular matrix $U$ as follows $M=LU$.
If $L$ has $1's$ on its diagonal, then it is called a Doolittle factorization. If
$U$ has $1's$ on its diagonal, then it is called a Crout factorization.
\end{definition}
When $M=LU$, we say that $M$ admits an $LU$-decomposition.\\
If $M=LU$ has a Doolittle factorization then we have the following picture:
\[
  \begin{bmatrix}
     m_{1,1} & m_{1,2} &  m_{1,3} & m_{1,4} & \cdots m_{1,n}\\
    m_{2,1} & m_{2,2} &  m_{2,3} & m_{2,4} & \cdots m_{2,n}\\
    m_{3,1} & m_{3,2} &  m_{3,3} & m_{3,4} & \cdots m_{3,n}\\
    m_{4,1} & m_{4,2} &  m_{4,3} & m_{4,4} & \cdots m_{4,n}\\
    \vdots\\
    m_{n,1} & m_{n,2} &  m_{n,3} & m_{n,4} & \cdots m_{n,n} \\

  \end{bmatrix}=(\begin{bmatrix}
    1       & 0  &  0   & 0&  \cdots & 0\\
    l_{2,1} &  1   &  0 & 0 & \cdots & 0\\
   l_{3,1} &  l_{3,2}   &  1 & 0& \cdots  & 0\\
    l_{4,1} &  l_{4,2}   &  l_{4,3} & 1 & \cdots  & 0\\
     \vdots\\
     l_{n,1} & l_{n,2} & l_{n,3} & l_{n,4} & \cdots & 1 \\
  \end{bmatrix},\begin{bmatrix}
     u_{1,1} & u_{1,2} &  u_{1,3} & u_{1,4} & \cdots u_{1,n}\\
    0 & u_{2,2} &  u_{2,3} & u_{2,4} & \cdots u_{2,n}\\
    0 & 0 &  u_{3,3} & u_{3,4} & \cdots u_{3,n}\\
    0 & 0 &  0 & u_{4,4} & \cdots u_{4,n}\\
    \vdots\\
     0 & 0 & 0 & 0 & \cdots & u_{n,n} \\

  \end{bmatrix})
\]

\section{$n$-matrix factorization of polynomials for $n=2$;$3$}

\subsection{2-matrix factorization of polynomials}

\textbf{Definition and some Examples}\\
Let $K[[x_{1},x_{2},\cdots, x_{m}]]$ be the power series ring in the indeterminates $x_{1},x_{2},\cdots, x_{m}$. We will sometimes write $K[[x]]$ instead of $K[[x_{1},x_{2},\cdots, x_{m}]]$ for ease of notation.\\
The notion of matrix factorization is defined in \cite{yoshino1998tensor} for nonzero non-invertible $f\in K[[x_{1},x_{2},\cdots, x_{m}]]$. We define it as in \cite{crisler2016matrix} slightly generalizing the one given in \cite{yoshino1998tensor} by including elements like $1\in K$ for convenience. Yoshino \cite{yoshino1998tensor} requires an element $f\in K[[x]]$ to be nonzero non-invertible because if $f=0$ then $K[[x]]/(f)=K[[x]]$ and if $f$ is a unit, then $K[[x]]/(f)=$K[[x]]/K[[x]]$=\{1\}$. But in this work, we will not bother about such restrictions because we will not deal with the homological methods used in \cite{yoshino1998tensor}.
\begin{definition}\cite{yoshino1998tensor}, \cite{crisler2016matrix} \cite{fomatati2019multiplicative} \label{defn matrix facto of polyn}   \\
An $m\times m$ \textbf{matrix factorization} of a polynomial $f\in \;R$ is a pair of $m$ $\times$ $m$ matrices $(P,Q)$ such that
$PQ=fI_{m}$, where $I_{m}$ is the $m \times m$ identity matrix and the coefficients of $P$ and of $Q$ are taken from $R$.
\end{definition}

 \begin{example} \label{exple: bad matrix facto of g}
 Let $l= xy+(x^{2}+yz)z$.
 We use the standard method (cf. \cite{fomatati2023refined}, \cite{fomatati2022tensor}, \cite{crisler2016matrix}) to find a matrix factorization of $l$ and quickly find:

  \[
  Q=(\begin{bmatrix}
    x & -(x^{2}+yz) \\
    z & y
  \end{bmatrix},\begin{bmatrix}
    y & x^{2}+yz \\
    -z & x
  \end{bmatrix})
\]

   Let $h= xy+x^{2}z+yz^{2}$. Observe that $l=h$. We can use the standard method to find a matrix factorization of $h$ with monomial entries. As we can see below, the matrix factors are nicer\footnote{"nicer" in the sense that we have matrices with monomial entries} than the ones obtained above. This approach was studied in \cite{fomatati2023refined}.\\
   First a matrix factorization of $xy+x^{2}z$ is
    \[
  (\begin{bmatrix}
    x & -x^{2} \\
    z & y
  \end{bmatrix},\begin{bmatrix}
    y & x^{2} \\
    -z & x
  \end{bmatrix})
\]
So, a matrix factorization of $h= xy+x^{2}z+yz^{2}$ is then:
 \[
  P=(\begin{bmatrix}
    x & -x^{2} & -y & 0\\
    z &   y    &  0 & -y\\
 z^{2} &  0    &  y &  x^{2}\\
     0 &  z^{2}&  -z & x
  \end{bmatrix},\begin{bmatrix}
     y & x^{2} &  y & 0\\
    -z &   x   &  0 & y\\
 -z^{2} &  0    &  x &  -x^{2}\\
     0 &  -z^{2}&  z & y
  \end{bmatrix})
\]
\end{example}

In this paper, we will not bother much about the type of matrix factors we will obtain, that is; whether they have monomial entries or not.
%
%
%
\\
We will simply be discussing how to factorize a polynomial using two or more matrices. We will refer to the type of factorizations of definition \ref{defn matrix facto of polyn} as 2-matrix factorizations because we have two matrix factors. This is the type that one easily finds in the literature (e.g. \cite{yoshino1998tensor}, \cite{crisler2016matrix}). We will generalize definition \ref{defn matrix facto of polyn} below (see definition \ref{defn: 3-matrix facto}).\\
In subsection \ref{example of 3-matrix facto}, we will show how to obtain a 3-matrix factorization of a polynomial $f$ from any of its 2-matrix factorizations. To that end, we will need to decompose one of the two matrix factors into two matrices using the Doolittle or the Crout matrix decomposition methods
that we will apply to matrices whose entries are not real numbers but polynomials. We recalled these methods in subsection \ref{Doolittle and Crout methods}.

\subsection{ 3-matrix factorizations and examples } \label{example of 3-matrix facto}
\begin{definition} \label{defn: 3-matrix facto}
  An $m\times m$ \textbf{$3$-matrix factorization} of a polynomial $f\in \;R$ is a triplet of $m$ $\times$ $m$ matrices $(A_{1},A_{2}, A_{3})$ such that
$A_{1} A_{2} A_{3}=fI_{m}$, where $I_{m}$ is the $m \times m$ identity matrix and the coefficients of each matrix $A_{i}$, $i\in \{1,2,3\}$, is taken from the field of fraction of $R$.
\end{definition}
In the following example, we find 3-matrix factors of a polynomial in $R$ which have entries in the field of fractions of $R$. That is why in the foregoing definition, we talk of the field of fraction or $R$ instead of just $R$ itself as we did in definition \ref{defn matrix facto of polyn}.

\begin{example}  \label{first exple of 3-matrix facto}

  Let $f=x^{2}+y^{2}$. \\A $2$-matrix factorization of $f$ is \[
  (\begin{bmatrix}
    x & -y \\
    y & x
  \end{bmatrix},\begin{bmatrix}
    x & y  \\
    -y & x
  \end{bmatrix})
=(x^{2}+y^{2})\begin{bmatrix}
    1 & 0  \\
    0 & 1
  \end{bmatrix}=fI_{2}.\]
  Thus, we have $fI_{2}=AB$ where $A=(\begin{bmatrix}
    x & -y \\
    y & x
  \end{bmatrix})$ and $B=(\begin{bmatrix}
    x & y  \\
    -y & x
  \end{bmatrix})$. \\

  Let us decompose $A$ as the product of a lower (L) and upper (U) triangular matrices, $A=LU$. We will use the Doolittle's approach (i.e., the main diagonal of $L$ shall be 1's).\\

  We have $\begin{bmatrix}
    x & -y \\
    y & x
  \end{bmatrix}=\begin{bmatrix}
    1   &   0 \\
    l_{21} & 1
  \end{bmatrix}\begin{bmatrix}
    U_{11} & U_{12}  \\
    0    & U_{22}
  \end{bmatrix}=\begin{bmatrix}
    U_{11} & U_{12}  \\
    l_{21}U_{11}    & l_{21}U_{12}+U_{22}
  \end{bmatrix}$

  Hence, $U_{11}=x$, $U_{12}=-y$. $l_{21}U_{11}=y \,\Rightarrow \,l_{21}=\frac{y}{x}$.\\
  $l_{21}U_{12}+U_{22}=x \, \Rightarrow \, U_{22}=x+\frac{y^{2}}{x}$.\\

  Hence, we obtain a $3$-matrix factorization of $f$: \\
  $\begin{bmatrix}
    1   &   0 \\
    \frac{y}{x} & 1
  \end{bmatrix}\begin{bmatrix}
    x & -y \\
    0    & x+\frac{y^{2}}{x}
  \end{bmatrix}\begin{bmatrix}
   x & -y  \\
    y    & x
  \end{bmatrix}=(x^{2}+y^{2})\begin{bmatrix}
    1   &   0 \\
    0 & 1
  \end{bmatrix}=fI_{2}$
\end{example}

\begin{remark}
Note that we could decompose $B$ instead of $A$ in order to obtain a $3$-matrix factorization of $f$. Also note that the Crout decomposition technique could be used in place of the Doolittle decomposition method.
\end{remark}
This remark is also valid for the following example.
\begin{example} \label{second exple of 3-matrix facto}
   Let $g=xyz+zx^{2}$. \\A $2$-matrix factorization of $g$ is \[
  (\begin{bmatrix}
    xy & -z \\
    x^{2} & z
  \end{bmatrix},\begin{bmatrix}
    z & z  \\
    -x^{2} & xy
  \end{bmatrix})
=(xyz+zx^{2})\begin{bmatrix}
    1 & 0  \\
    0 & 1
  \end{bmatrix}=gI_{2}.\]

  Thus, we have $gI_{2}=PQ$ where $P=(\begin{bmatrix}
      xy & -z \\
    x^{2} & z
  \end{bmatrix})$ and $Q=(\begin{bmatrix}
    z & z  \\
    -x^{2} & xy
  \end{bmatrix})$.\\

   Let us decompose $P$ as the product of a lower (L) and upper (U) triangular matrices, $P=LU$. We will again use the Doolittle's approach.\\

  We have $\begin{bmatrix}
   xy & -z \\
    x^{2} & z
  \end{bmatrix}=\begin{bmatrix}
    1   &   0 \\
    l_{21} & 1
  \end{bmatrix}\begin{bmatrix}
    U_{11} & U_{12}  \\
    0    & U_{22}
  \end{bmatrix}=\begin{bmatrix}
    U_{11} & U_{12}  \\
    l_{21}U_{11}    & l_{21}U_{12}+U_{22}
  \end{bmatrix}$

  Hence, $U_{11}=xy$, $U_{12}=-z$. $l_{21}U_{11}=x^{2} \,\Rightarrow \,l_{21}=\frac{x}{y}$.\\
  $l_{21}U_{12}+U_{22}=z \, \Rightarrow \, U_{22}=z+\frac{zx}{y}$.\\

  Hence, we obtain a $3$-matrix factorization of $g$: \\
  $\begin{bmatrix}
    1   &   0 \\
    \frac{x}{y} & 1
  \end{bmatrix}\begin{bmatrix}
    xy & -z \\
    0    & z+\frac{zx}{y}
  \end{bmatrix}\begin{bmatrix}
  z & z  \\
    -x^{2} & xy
  \end{bmatrix}=(xyz+zx^{2})\begin{bmatrix}
    1   &   0 \\
    0 & 1
  \end{bmatrix}=gI_{2}$
 \\\\
\end{example}

\section{The category of 3-matrix factorizations of $f\in R$} \label{subsec category of 3-matrix facto of f}

In this section, we construct the category of $3$-matrix factorizations of $f\in R$. In the sequel, $K(x_{1},\cdots,x_{n})$ denotes the fraction field of $K[x_{1},\cdots,x_{n}]$. \\
The category of $3$-matrix factorizations of a polynomial $f\in R=K[x]:=K[x_{1},\cdots,x_{n}]$ which we denote by $MF(R,f)_{3}$ or $MF_{R}(f)_{3}$, or $MF(f)_{3}$ (when there is no risk of confusion) is defined as follows:\\
$\bullet$ The objects are the $3$-matrix factorizations of $f$.\\
$\bullet$ Given two $3$-matrix factorizations of $f$; $(\phi_{1},\psi_{1},\theta_{1})$ and $(\phi_{2},\psi_{2},\theta_{2})$ respectively of sizes $n_{1}$ and $n_{2}$, a morphism from $(\phi_{1},\psi_{1},\theta_{1})$ to $(\phi_{2},\psi_{2},\theta_{2})$ is a triplet of matrices $(\alpha,\beta,\delta)$ each of size $n_{2}\times n_{1}$ which makes the following diagram commute:

$$\xymatrix@ R=0.6in @ C=.75in{K(x)^{n_{1}} \ar[r]^{\theta_{1}} \ar[d]_{\alpha} &K(x)^{n_{1}} \ar[r]^{\psi_{1}} \ar[d]_{\delta} &
K(x)^{n_{1}} \ar[d]^{\beta} \ar[r]^{\phi_{1}} & K(x)^{n_{1}}\ar[d]^{\alpha \;\;\;\;\;\;\;\;\;\;}\\
K(x)^{n_{2}} \ar[r]^{\theta_{2}} &K(x)^{n_{2}} \ar[r]^{\psi_{2}} & K(x)^{n_{2}}\ar[r]^{\phi_{2}} & K(x)^{n_{2}}}$$
%
That is,
$$\begin{cases}
 \alpha\phi_{1}=\phi_{2}\beta  \\
 \psi_{2}\delta= \beta\psi_{1} \\
 \delta\theta_{1}=\theta_{2}\alpha
\end{cases}$$
$\bullet$ Given three $3$-matrix factorizations of $f$: $(\phi_{1},\psi_{1},\theta_{1})$, $(\phi_{2},\psi_{2},\theta_{2})$ and $(\phi_{3},\psi_{3},\theta_{3})$ respectively of sizes $n_{1},n_{2}$ and $n_{3}$, the composition

$$\xymatrix{(\alpha_{2},\beta_{2},\delta_{2})\circ (\alpha_{1},\beta_{1},\delta_{1}): (\phi_{1},\psi_{1},\theta_{1}) \ar[r]^{\,\,\,\,\,\,\,\,\,\,\,\,\,\,\,\,\,\,\,\,\,\,\,\,\,\,\,\,\,\,\,\,\,\,\,\,\,\,\,\,\,\,\,\,\,\,(\alpha_{1},\beta_{1},\delta_{1})} & (\phi_{2},\psi_{2},\theta_{2}) \ar[r] ^{(\alpha_{2},\beta_{2},\delta_{2})} & (\phi_{3},\psi_{3},\theta_{3})}$$
is the triplet of matrices $(\alpha_{2}\alpha_{1},\beta_{2}\beta_{1},\delta_{2}\delta_{1})$ such that the following diagram commutes:
$\xymatrix@ R=0.6in @ C=.75in{K(x)^{n_{1}} \ar[r]^{\theta_{1}} \ar[d]_{\alpha_{1}} &K(x)^{n_{1}} \ar[r]^{\psi_{1}} \ar[d]_{\delta_{1}} &
K(x)^{n_{1}} \ar[d]^{\beta_{1}} \ar[r]^{\phi_{1}} & K(x)^{n_{1}}\ar[d]^{\alpha_{1}}\\
K(x)^{n_{2}} \ar[r]^{\theta_{2}} & K(x)^{n_{2}} \ar[r]^{\psi_{2}} & K(x)^{n_{2}}\ar[r]^{\phi_{2}} & K(x)^{n_{2}}
}\\
\xymatrix@ R=0.6in @ C=.75in{ \ar[d]_{\alpha_{2}} & \ar[d]_{\delta_{2}} &
 \ar[d]^{\beta_{2}} & \ar[d]^{\alpha_{2}}\\
K(x)^{n_{3}} \ar[r]^{\theta_{3}} & K(x)^{n_{3}} \ar[r]^{\psi_{3}} & K(x)^{n_{3}}\ar[r]^{\phi_{3}} & K(x)^{n_{3}}
}$
%

That is,
$$\begin{cases}
 (\alpha_{2}\alpha_{1})\phi_{1}=\phi_{3}(\beta_{2}\beta_{1})  \\
 \psi_{3}(\delta_{2}\delta_{1})= (\beta_{2}\beta_{1})\psi_{1} \\
 (\delta_{2}\delta_{1})\theta_{1}= \theta_{3}(\alpha_{2}\alpha_{1})
\end{cases}$$
$\bullet$ It is easy to see that \textit{associativity of composition of maps of $3$-matrix factorizations of $f$} is a consequence of the fact that matrix multiplication is associative.\\
$\bullet$ For any $n\times n$ $3$-matrix factorization $(\phi,\psi,\theta)$ of $f$, there is a map $1_{(\phi,\psi,\theta)}: (\phi,\psi,\theta) \rightarrow (\phi,\psi,\theta)$ which is actually the triplet of identity $n\times n$ matrices $(I_{n},I_{n},I_{n})$.
\\
$\bullet$ It is also clear that composing any map of $3$-matrix factorizations of $f$ with $1_{(\phi,\psi,\theta)}$ from the left or the right (whenever the composition is possible) leaves the given map unchanged. This ends the definition of the category of $3$-matrix factorizations of $f\in R=K[x]$.

\section{The multiplicative tensor product of 3-matrix factorizations}

In this section, we define of the multiplicative tensor product of $3$-matrix factorizations denoted $\overline{\otimes}_{3}$ and we prove that it is a bifunctorial operation.  The functoriality of this new operation will be proved in subsection \ref{sec: functoriality of 3-matrix tens prod}. Finally, we give some properties of this operation.

\subsection{Definition and examples} \label{subsec: defn and properties}

\begin{definition}   \label{defn of the mult tens prodt of 3-matr facto}
Let $X=(\phi,\psi,\theta)$ be a matrix factorization of $f\in K[x]$ of size $n$ and let $X'=(\phi',\psi',\theta')$ be a matrix factorization of $g\in K[y]$ of size $m$. Thus, $\phi,\psi,\phi',\psi',\theta,\,and\, \theta'$ can be considered as matrices over $L=K(x,y)$. The \textbf{multiplicative tensor product of 3-matrix factorizations} $X\overline{\otimes}_{3} X'$ is given by \\\\
$X\overline{\otimes}_{3} X'$ = $(\phi,\psi,\theta)\overline{\otimes}_{3} (\phi',\psi',\theta')=([\phi\otimes\phi'], [\psi\otimes\psi'], [\theta\otimes\theta'])$
\\\\
where each component is an endomorphism on $L^{n}\otimes_{L} L^{m}$.
\end{definition}

\begin{example}
  Consider $g=xyz+zx^{2}$ and $f=x^{2}+y^{2}$ the polynomials of examples \ref{first exple of 3-matrix facto} and \ref{second exple of 3-matrix facto}.\\
  $X=(\begin{bmatrix}
    1   &   0 \\
    \frac{y}{x} & 1
  \end{bmatrix},\begin{bmatrix}
    x & -y \\
    0    & x+\frac{y^{2}}{x}
  \end{bmatrix},\begin{bmatrix}
   x & y  \\
    -y  & x
  \end{bmatrix})$ is a $3$-matrix factorization of $f$.\\
  $X'=(\begin{bmatrix}
    1   &   0 \\
    \frac{x}{y} & 1
  \end{bmatrix},\begin{bmatrix}
    xy & -z \\
    0    & z+\frac{zx}{y}
  \end{bmatrix},\begin{bmatrix}
  z & z  \\
    -x^{2} & xy
  \end{bmatrix})$ is a $3$-matrix factorization of $g$.\\\\

  $X\overline{\otimes}_{3} X'$ = $(\begin{bmatrix}
    1   &   0 \\
    \frac{y}{x} & 1
  \end{bmatrix},\begin{bmatrix}
    x & -y \\
    0    & x+\frac{y^{2}}{x}
  \end{bmatrix},\begin{bmatrix}
   x & y  \\
    -y  & x
  \end{bmatrix})\overline{\otimes}_{3} (\begin{bmatrix}
    1   &   0 \\
    \frac{x}{y} & 1
  \end{bmatrix},\begin{bmatrix}
    xy & -z \\
    0    & z+\frac{zx}{y}
  \end{bmatrix},\begin{bmatrix}
  z & z  \\
    -x^{2} & xy
  \end{bmatrix})\\=(\begin{bmatrix}
    1   &   0 \\
    \frac{y}{x} & 1
  \end{bmatrix}\otimes\begin{bmatrix}
    1   &   0 \\
    \frac{x}{y} & 1
  \end{bmatrix}, \begin{bmatrix}
    x & -y \\
    0    & x+\frac{y^{2}}{x}
  \end{bmatrix}\otimes\begin{bmatrix}
    xy & -z \\
    0    & z+\frac{zx}{y}
  \end{bmatrix}, \begin{bmatrix}
   x &  y  \\
    -y  & x
  \end{bmatrix}\otimes\begin{bmatrix}
  z & z  \\
    -x^{2} & xy
  \end{bmatrix})\\=(\begin{bmatrix}
    1   &   0       &0     &0  \\
    \frac{x}{y} & 1  &0     &0 \\
    \frac{y}{x}  & 0     &1    &0 \\
    1      &  \frac{y}{x}  &\frac{x}{y}  &1
  \end{bmatrix}, \begin{bmatrix}
    x^{2}y   &   -xz       &-xy^{2}     &yz  \\
    0       & xz+\frac{zx^{2}}{y}  &0     &-zy-zx \\
    0       & 0     & x^{2}y+y^{3}    & -zx-\frac{zy^{2}}{x} \\
    0      &  0  & 0  &   xz+\frac{zx^{2}}{y}+\frac{y^{2}z}{x}+zy
  \end{bmatrix}, \begin{bmatrix}
    xz   &   xz       & yz     & yz  \\
    -x^{3} &   x^{2}y  & -x^{2}y     & xy^{2} \\
    -yz  &  -yz    & xz    & xz \\
    yx^{2}   &  -xy^{2}  & -x^{3}  & x^{2}y
  \end{bmatrix})$ is a $3$-matrix factorization of $fg$.\\\\
  In fact,\\\\ $\begin{bmatrix}
    1   &   0       &0     &0  \\
    \frac{x}{y} & 1  &0     &0 \\
    \frac{y}{x}  & 0     &1    &0 \\
    1      &  \frac{y}{x}  &\frac{x}{y}  &1
  \end{bmatrix} \begin{bmatrix}
    x^{2}y   &   -xz       &-xy^{2}     &yz  \\
    0       & xz+\frac{zx^{2}}{y}  &0     &-zy-zx \\
    0       & 0     & x^{2}y+y^{3}    & -zx-\frac{zy^{2}}{x} \\
    0      &  0  & 0  &   xz+\frac{zx^{2}}{y}+\frac{y^{2}z}{x}+zy
  \end{bmatrix} \begin{bmatrix}
    xz   &   xz       & yz     & yz  \\
    -x^{3} &   x^{2}y  & -x^{2}y     & xy^{2} \\
    -yz  &  -yz    & xz    & xz \\
    yx^{2}   &  -xy^{2}  & -x^{3}  & x^{2}y
  \end{bmatrix}=\\ \begin{bmatrix}
    x^{2}y   &   -xz       &-xy^{2}     &yz  \\
    x^{3}  & \frac{-x^{2}z}{y}+xz+\frac{zx^{2}}{y}  & -x^{2}y     &-zy \\
    xy^{2}       & -yz    & x^{2}y    &\frac{y^{2}z}{x}-zx-\frac{zy^{2}}{x} \\
    x^{2}y      &  yz  & x^{3}  & xz
  \end{bmatrix} \begin{bmatrix}
    xz   &   xz       & yz     & yz  \\
    -x^{3} &   x^{2}y  & -x^{2}y     & xy^{2} \\
    -yz  &  -yz    & xz    & xz \\
    yx^{2}   &  -xy^{2}  & -x^{3}  & x^{2}y
  \end{bmatrix}=fgI_{4}$

\end{example}

\subsection{Funtoriality of $\overline{\otimes}_{3}$} \label{sec: functoriality of 3-matrix tens prod}
This subsection is entirely devoted to the discussion of the bifunctoriality of $\overline{\otimes}_{3}$. \\

  \textbf{Setting the stage:}
Let $X_{f}=(\phi,\psi,\theta)$, $X'_{f}=(\phi',\psi',\theta')$ and $X_{f}''=(\phi'',\psi'',\theta'')$ be objects of $MF(K[x],f)_{3}$ respectively of sizes $n, n'$ and $n''$. Let $X_{g}=(\sigma,\rho,\zeta )$, $X_{g}'=(\sigma',\rho',\zeta')$ and $X_{g}''=(\sigma'',\rho'',\zeta'')$ be objects of $MF(K[y],g)_{3}$ respectively of sizes $m, m'$ and $m''$.

\begin{definition}\label{defn zeta is a bifuntor}
For morphisms $\Phi_{f}=(\alpha_{f}, \beta_{f},\delta_{f}): X_{f}=(\phi,\psi,\theta) \rightarrow X_{f}'=(\phi',\psi',\theta')$  and $\Phi_{g}=(\alpha_{g}, \beta_{g},\delta_{g}): X_{g}=(\sigma,\rho,\zeta) \rightarrow X_{g}'=(\sigma',\rho',\zeta')$ respectively in $MF(K[x],f)_{3}$ and $MF(K[y],g)_{3}$,\\ we define $\Phi_{f}\overline{\otimes}_{3} \Phi_{g}: X_{f}\overline{\otimes} _{3} X_{g} =(\phi,\psi,\theta)\overline{\otimes}_{3} (\sigma,\rho,\zeta)\rightarrow X_{f}'\overline{\otimes}_{3} X_{g}' =(\phi',\psi',\theta')\overline{\otimes}_{3} (\sigma',\rho',\zeta')$
by
$$([\alpha_{f}\otimes \alpha_{g}],[\beta_{f}\otimes \beta_{g}],[\delta_{f}\otimes \delta_{g}])
 $$
\end{definition}

\begin{lemma}\label{zeta is a morphism in both arguments}

  $\Phi_{f}\overline{\otimes}_{3} \Phi_{g}$: $X_{f}\overline{\otimes}_{3}X_{g} =(\phi,\psi,\theta)\overline{\otimes}_{3} (\sigma,\rho,\zeta)\rightarrow X_{f}'\overline{\otimes}_{3}X_{g}' =(\phi',\psi',\theta')\overline{\otimes}_{3} (\sigma',\rho',\zeta')$ is a morphism in $MF(K[x,y],fg)_{3}$.

\end{lemma}
\begin{proof}
 We need to show that the following diagram commutes:


$$\xymatrix@ R=0.8in @ C=1.10in{K(x,y)^{nm} \ar[r]^{[\theta\otimes\zeta ]} \ar[d]_{[ \alpha_{f} \otimes \alpha_{g}]} & K(x,y)^{nm} \ar[r]^{[\psi\otimes \rho]} \ar[d]_{[\delta_{f} \otimes \delta_{g}]} &
K(x,y)^{nm} \ar[d]^{[\beta_{f} \otimes \beta_{g}]} \ar[r]^{[\phi\otimes\sigma]} & K(x,y)^{nm}\ar[d]^{[\alpha_{f} \otimes \alpha_{g}]}\\
K(x,y)^{n'm'}\ar[r]_{[\theta'\otimes\zeta']} & K(x,y)^{n'm'} \ar[r]_{[\psi'\otimes\rho']} & K(x,y)^{n'm'}\ar[r]_{[\phi'\otimes\sigma']} & K(x,y)^{n'm'}}$$

%

 viz. all the three squares in the foregoing diagram commute.\\
 $\bullet$ The commutativity of each of these squares from the right to the left is expressed by the following equalities:\\
 $ [\alpha_{f} \otimes \alpha_{g} ][\phi\otimes\sigma]=[\phi'\otimes\sigma'][\beta_{f} \otimes \beta_{g}]$
\\
 $[\beta_{f} \otimes \beta_{g}][\psi\otimes\rho]=[\psi'\otimes\rho'][\delta_{f}\otimes \delta_{g}]$ and\\
$[\delta_{f} \otimes \delta_{g}][\theta\otimes\zeta]=[\theta'\otimes\zeta'][\alpha_{f} \otimes \alpha_{g}]$

i.e., all we need to show is the set of equalities:
\\
$$\begin{cases}
\alpha_{f}\phi\otimes \alpha_{g}\sigma=\phi'\beta_{f}\otimes\sigma'\beta_{g} \cdots (1)\\
\beta_{f}\psi \otimes \beta_{g}\rho= \psi'\delta_{f}\otimes\rho'\delta_{g} \cdots (2)\\
\delta_{f}\theta \otimes \delta_{g}\zeta= \theta'\alpha_{f}\otimes\zeta'\alpha_{g} \cdots (3)
\end{cases}$$

Now by hypothesis, $\Phi_{f}=(\alpha_{f},\beta_{f},\delta_{f}): X_{f}=(\phi,\psi,\theta) \rightarrow X_{f}'=(\phi',\psi',\theta')$ and $\Phi_{g}=(\alpha_{g},\beta_{g},\delta_{g}): X_{g}=(\sigma,\rho,\zeta) \rightarrow X_{g}'=(\sigma',\rho',\zeta')$ are morphisms, meaning that the following diagrams commute
$$\xymatrix@ R=0.6in @ C=.75in{K(x)^{n} \ar[r]^{\theta} \ar[d]_{\alpha_{f}} &K(x)^{n} \ar[r]^{\psi} \ar[d]_{\delta_{f}} &
K(x)^{n} \ar[d]^{\beta_{f}} \ar[r]^{\phi} & K(x)^{n}\ar[d]^{\alpha_{f} }\\
K(x)^{n'} \ar[r]^{\theta'} &K(x)^{n'} \ar[r]^{\psi'} & K(x)^{n'}\ar[r]^{\phi'} & K(x)^{n'}}$$
%

and \\
$$\xymatrix@ R=0.6in @ C=.75in{K(y)^{m} \ar[r]^{\zeta} \ar[d]_{\alpha_{g}} &K(y)^{m} \ar[r]^{\rho} \ar[d]_{\delta_{g}} &
K(y)^{m} \ar[d]^{\beta_{g}} \ar[r]^{\sigma} & K(y)^{m}\ar[d]^{\alpha_{g} }\\
K(y)^{m'} \ar[r]^{\zeta'} &K(y)^{m'} \ar[r]^{\rho'} & K(y)^{m'}\ar[r]^{\sigma'} & K(y)^{m'}}$$
  That is,
$$\begin{cases}
 \alpha_{f}\phi=\phi'\beta_{f} \cdots (i) \\
 \psi'\delta_{f}= \beta_{f}\psi \cdots (ii) \\
 \theta'\alpha_{f}= \delta_{f}\theta \cdots (iii)
\end{cases}$$
and \\
$$\begin{cases}
 \alpha_{g}\sigma=\sigma'\beta_{g} \cdots (i') \\
 \rho'\delta_{g}= \beta_{g}\rho \cdots (ii') \\
 \zeta'\alpha_{g}= \delta_{g}\zeta \cdots (iii')
\end{cases}$$
Now considering $(i)$ and $(i')$, we immediately see that equality $(1)$ holds. Similarly, $(ii)$ and $(ii')$ yield $(2)$. Finally, $(iii)$ and $(iii')$ yield $(3)$\\

So, $\Phi_{f} \overline{\otimes}_{3} \Phi_{g}$ is a morphism in $MF(K[x,y], fg)$.

\end{proof}
We can now state the following result.
\begin{theorem} \label{zeta is a bifunctor}
\begin{enumerate}
\item Let $X$ be a matrix factorization of $f\in K[x]$ of size $n$ and let $Y$ be a matrix factorization of $g\in K[y]$ of size $m$. Then, there is a tensor product $\overline{\otimes}_{3}$ of $3$-matrix factorizations which produces a $3$-matrix factorization $X \overline{\otimes}_{3} Y$ of the product $fg\in K[x_{1},...,x_{r},y_{1},...,y_{s}]$ which is of size $nm$. $\overline{\otimes}_{3}$ is called the multiplicative tensor product of $3$-matrix factorizations.
\item The multiplicative tensor product of $3$-matrix factorizations $(-) \overline{\otimes}_{3} (-):MF(K[x],f)\times MF(K[y],g)\rightarrow MF(K[x,y],fg)$ is a bifunctor.
\end{enumerate}
\end{theorem}
\begin{proof}
\begin{enumerate}
  \item This is exactly what we proved above in subsection \ref{subsec: defn and properties}.
  \item We show that $\overline{\otimes}_{3}$ is a bifunctor.\\
In order to ease our computations, let's write $F=(-)\overline{\otimes}_{3} (-)$. We show that $F$ is a bifunctor.
\newpage
 We have:

$\;\; (-)\overline{\otimes}_{3} (-):\;\;\;\;\;\;\;\;\;\;\;\;\; MF(f)\times MF(g)\,\,\,\,\,\,\,\,\,\,\longrightarrow \,\,\,\,\,\,\,\,\,\,\,\,\,\,\,\,MF(fg)$
$$\xymatrix @ R=0.4in @ C=.33in
{&(X_{f}\;\;\;\;\;\;{,} \ar[d]_{\Phi_{f}}\ar @{}[dr]& X_{g}) \ar[rrr]^-{} \ar[d]_{\Phi_{g}}
\ar @{}[dr] &&& X_{f}\overline{\otimes}_{3} X_{g} \ar[d]^{\Phi_{f}\overline{\otimes}_{3} \Phi_{g}:= (\alpha,\beta,\delta)}\\
&(X_{f}'\;\;\;\;\;\; {,}& X_{g}') &\ar[r]&& X_{f}'\overline{\otimes}_{3} X_{g}'}$$
$$\xymatrix @ R=0.4in @ C=.3in
{&\ar[d]_{\Phi_{f}'}\ar @{}[dr]& \ar[d]_{\Phi_{g}'}
\ar @{}[dr] &&&\ar[d]^{\Phi_{f}'\overline{\otimes}_{3} \Phi_{g}':= (\alpha',\beta',\delta')}\\
&(X_{f}'' \;\;\;\;\;\;{,}& X_{g}'') &\ar[r]&& X_{f}''\overline{\otimes}_{3} X_{g}''}$$

We showed in lemma \ref{zeta is a morphism in both arguments} that $\Phi_{f}\overline{\otimes}_{3} \Phi_{g}:= (\alpha,\beta,\delta)$ is a morphism in
$MF(K[x,y],fg)$, where
 $$(\alpha,\beta,\delta)=([\alpha_{f}\otimes \alpha_{g}],[\beta_{f}\otimes \beta_{g}],[\delta_{f}\otimes \delta_{g}])$$
Similarly, if $\Phi_{f}':=(\alpha_{f}',\beta_{f}',\delta_{f}')$ and $\Phi_{g}':=(\alpha_{g}',\beta_{g}',\delta_{g}')$ then $\Phi_{f}' \overline{\otimes}_{3} \Phi_{g}'=(\alpha',\beta',\delta')$ where $$(\alpha',\beta',\delta')=([\alpha_{f}'\otimes \alpha_{g}'],[\beta_{f}'\otimes \beta_{g}'],[\delta_{f}'\otimes \delta_{g}'])$$
 It now remains to show the composition and the identity axioms.\\\\
\textit{Identity Axiom}: \\
We show that $F(id_{(X_{f},X_{g})})=id_{F(X_{f},X_{g})}$.\\
Now, $F(id_{(X_{f},X_{g})})=F(id_{X_{f}},id_{X_{g}}):=id_{X_{f}}\overline{\otimes}_{3} id_{X_{g}}: X_{f} \overline{\otimes}_{3}X_{g} \rightarrow X_{f} \overline{\otimes}_{3}X_{g}$.\\
And by definition \ref{defn zeta is a bifuntor}, $id_{X_{f}}\overline{\otimes}_{3} id_{X_{g}}$ is the triplet of matrices \\
 $([I_{n}\otimes I_{m}],[I_{n}\otimes I_{m}],[I_{n}\otimes I_{m}])\,\,\,\,\,\,\,\,\,\,\,\dag$ \\

Next, we compute $id_{F(X_{f},X_{g})}=id_{X_{f} \overline{\otimes}_{3}X_{g}}:X_{f} \overline{\otimes}_{3}X_{g} \rightarrow X_{f} \overline{\otimes}_{3}X_{g}$.\\
By definition of a morphism in the category $MF(fg)$, we know that \\

$id_{X_{f} \overline{\otimes}_{3}X_{g}}:=([I_{nm}],[I_{nm}],[I_{nm}])\,\,\,\,\,\,\,\,\,\,\,\dag\dag
$\\

Since $I_{n}\otimes I_{m}= I_{nm}$, we see that $\dag$ and $\dag\dag$ are the same, therefore $F(id_{(X_{f},X_{g})})=id_{F(X_{f},X_{g})}$ as desired.\\\\
\textit{Composition Axiom}:\\
Consider the situation:
$$\xymatrix@R=4in{X_{f}\ar[r]^{\Phi_{f}} & X_{f}'\ar[r]^{\Phi_{f}'}& X_{f}''}$$ $$\xymatrix@R=4in{X_{g}\ar[r]^{\Phi_{g}} & X_{g}'\ar[r]^{\Phi_{g}'}& X_{g}''}$$
$$\xymatrix@R=4in{X_{f} \overline{\otimes}_{3}X_{g}\ar[r]^{F(\Phi_{f},\Phi_{g})} & X_{f}'\overline{\otimes}_{3}X_{g}'\ar[r]^{F(\Phi_{f}',\Phi_{g}')}& X_{f}''\overline{\otimes}_{3} X_{g}''}$$

We need to show $F(\Phi_{f}'\circ \Phi_{f},\Phi_{g}'\circ \Phi_{g})=F(\Phi_{f}',\Phi_{g}')\circ F(\Phi_{f},\Phi_{g})$. \\
Now, $\Phi_{f}'\circ \Phi_{f}= (\alpha_{f}'\alpha_{f},\beta_{f}'\beta_{f})$ and $\Phi_{g}'\circ \Phi_{g}=(\alpha_{g}'\alpha_{g},\beta_{g}'\beta_{g})$.\\\\
Thanks to definition \ref{defn zeta is a bifuntor}, we obtain:
$$(\Phi_{f}'\circ \Phi_{f})\overline{\otimes}_{3}(\Phi_{g}'\circ \Phi_{g})=([\alpha_{f}'\alpha_{f}\otimes \alpha_{g}'\alpha_{g}],[\beta_{f}'\beta_{f}\otimes \beta_{g}'\beta_{g}],[\delta_{f}'\delta_{f}\otimes \delta_{g}'\delta_{g}])\,\,\,\,\,\,\,\,\,\,\,\, \ddag'$$

Next, \\\\
$(\Phi_{f}'\overline{\otimes}_{3}\Phi_{g}')\circ (\Phi_{f}\overline{\otimes}_{3} \Phi_{g})\\
=([\alpha_{f}'\otimes \alpha_{g}'],[\beta_{f}'\otimes \beta_{g}'],[\delta_{f}'\otimes \delta_{g}'])
\circ
([\alpha_{f}\otimes \alpha_{g}],[\beta_{f}\otimes \beta_{g}],[\delta_{f}\otimes \delta_{g}])
\\
=([\alpha_{f}'\alpha_{f}\otimes \alpha_{g}'\alpha_{g}],[\beta_{f}'\beta_{f}\otimes \beta_{g}'\beta_{g}],[\delta_{f}'\delta_{f}\otimes \delta_{g}'\delta_{g}])
\,\,\,\,\,\,\,\,\,\,\,\, \ddag\ddag'$

From $\ddag'$ and $\ddag\ddag'$, we see that $F(\Phi_{f}'\circ \Phi_{f},\Phi_{g}'\circ \Phi_{g})=F(\Phi_{f}',\Phi_{g}')\circ F(\Phi_{f},\Phi_{g})$.
Thus, $(-)\overline{\otimes}_{3}(-)$ is a bifunctor.

\end{enumerate}
\end{proof}

\subsection{Properties of $\overline{\otimes}_{3}$: Associativity, commutativity and distributivity}

In this subsection, we give some properties of the multiplicative tensor product of 3-matrix factorizations. We prove that $\overline{\otimes}_{3}$ is associative, commutative and distributive.
\\
We denote by $X_{1}=(\phi_{1},\psi_{1},\theta_{1})$ (resp. $X_{2}=(\phi_{2},\psi_{2},\theta_{2})$) an $(n_{1}\times n_{1})$ (resp. $(n_{2}\times n_{2})$) $3$-matrix factorization of $f\in K[x]$. We also let $X'=(\phi',\psi',\theta')$ (resp. $X''=(\phi'',\psi'',\theta'')$) denotes a $(p\times p)$ (resp. $(m\times m)$) $3$-matrix factorization of $g\in K[y]$ (resp. of $h\in K[z]:= K[z_{1},\cdots,z_{l}]$).
$X=(\phi,\psi,\theta)$ will also be an $r\times r$ $3$-matrix factorization of $f\in K[x]$.

\begin{proposition} (Associativity) \label{Assoc of new tensor product}\\
There is an identity:\\
  $(X\overline{\otimes}_{3}X')\overline{\otimes}_{3}X''=X\overline{\otimes}_{3}(X'\overline{\otimes}_{3}X'')$ in $MF(fgh)$.

\end{proposition}
\begin{proof}
The desired identity follows from the fact that the standard tensor product for matrices is associative.
\end{proof}
To prove the commutativity of $\overline{\otimes}_{3}$, recall (cf. section 3.1 \cite{henderson1981vec}) that given two matrices $C$ and $D$, the tensor products $C \otimes D$ and $D \otimes C$ are \textbf{permutation equivalent}. That is, there exist permutation matrices $P$ and $Q$ (so called commutation matrices) such that:
$C \otimes D = P (D \otimes C) Q$. If $C$ and $D$ are square matrices, then $C \otimes D$ and $D\otimes C$ are even \textbf{permutation similar}, meaning we can take $P = Q^{T}$.\\
To be more precise \cite{henderson1981vec}, if $C$ is a $p\times q$ matrix and $D$ is an $r\times s$ matrix, then $$D \otimes C = S_{p,r} (C \otimes D) S_{q,s}^{T}$$
where, $$ S_{m,n}=\sum_{i=1}^{m}( e_{i}^{T}\otimes I_{n}\otimes e_{i}) = \sum_{j=1}^{n}( e_{j}\otimes I_{m}\otimes e_{j}^{T})$$
$I_{n}$ is the $n\times n$ identity matrix and $e_{i}$ is the $i^{th}$ unit vector. $ S_{m,n}$ is the \textbf{perfect shuffle} permutation matrix.
\\
The commutativity of $\overline{\otimes}_{3}$ is up to isomorphism. This isomorphism comes from the permutation similarity
of the matrices involved.
\begin{proposition} (commutativity) \label{prop commutativity of new tensor product}\\
For $3$-matrix factorizations $X\in MF(f)$ and $X'\in MF(g)$, there is a natural isomorphism

   $X\overline{\otimes}_{3} X'\cong X' \overline{\otimes}_{3} X \,in\,MF(fg).$

\end{proposition}
\begin{proof}
We know that $X\overline{\otimes}_{3} X' =([\phi\otimes \phi'],[\psi\otimes \psi'],[\theta\otimes \theta'])\,
and\,X'\overline{\otimes}_{3} X=([\phi'\otimes \phi],[\psi'\otimes \psi],[\theta'\otimes \theta])$.
The desired isomorphism then follows from the fact that: $\phi\otimes \phi'$ (respectively $\psi\otimes \psi'$) and $\phi'\otimes \phi$ (respectively $\psi'\otimes \psi$) are permutation similar, and also that $\theta\otimes \theta'$ and $\theta'\otimes \theta$ are permutation equivalent.
\end{proof}
\begin{proposition} (Distributivity)\\
If $X_{1}$ and $X_{2}$ are $3$-matrix factorizations (of $f\in K[x]$) of the same size, then there are identities
\begin{enumerate}
  \item $(X_{1}\oplus X_{2})\overline{\otimes}_{3} X'=(X_{1}\overline{\otimes}_{3} X')\oplus (X_{2}\overline{\otimes}_{3} X').$
  \item $ X' \overline{\otimes}_{3}(X_{1}\oplus X_{2})=(X'\overline{\otimes}_{3}X_{1})\oplus (X'\overline{\otimes}_{3}X_{2}).$
  \end{enumerate}

\end{proposition}
\begin{proof}
\begin{enumerate}
  \item
$$(X_{1}\overline{\otimes}_{3} X')\oplus (X_{2}\overline{\otimes}_{3} X')$$
$\,\,\,\,\,\,\,\,\,\,\,\,\,\,\,\,\,\,\,\,\,\,\,\,\,\,\,\,\,\,\,\,\,\,\,\,\,\,=([\phi_{1}\otimes \phi'],[\psi_{1}\otimes \psi'],[\theta_{1}\otimes \theta'])\oplus ([\phi_{2}\otimes \phi'],[\psi_{2}\otimes \psi'],[\theta_{2}\otimes \theta'])$\\
\[=(\begin{bmatrix}
  \phi_{1}\otimes \phi' &   0      \\
  0   & \phi_{2}\otimes \phi'
\end{bmatrix}, \begin{bmatrix}
  \psi_{1}\otimes \psi' &   0     \\
  0   & \psi_{2}\otimes \psi'
\end{bmatrix}, \begin{bmatrix}
  \theta_{1}\otimes \theta' &   0     \\
  0   & \theta_{2}\otimes \theta'
\end{bmatrix})\,\,\,\,\,\,\,\,\cdots (\sharp)\]\\
Next, $$(X_{1}\oplus X_{2})\overline{\otimes}_{3} X'=((\phi_{1},\psi_{1},\theta_{1})\oplus (\phi_{2},\psi_{2},\theta_{2}))\overline{\otimes}_{3} (\phi',\psi',\theta')$$
\[=(\begin{bmatrix}
  \phi_{1} &          0      \\
    0               &   \phi_{2}
\end{bmatrix}, \begin{bmatrix}
  \psi_{1}   &          0      \\
    0               &   \psi_{2}
\end{bmatrix}, \begin{bmatrix}
  \theta_{1}   &          0      \\
    0               &   \theta_{2}
\end{bmatrix})\overline{\otimes}(\phi',\psi',\theta')\]\\
\[=(
  \begin{bmatrix}
  \phi_{1} &          0      \\
    0               &   \phi_{2}
\end{bmatrix}\otimes\phi' , \begin{bmatrix}

  \psi_{1}   &          0      \\
    0               &   \psi_{2}
\end{bmatrix}\otimes\psi',\begin{bmatrix}

  \theta_{1}   &          0      \\
    0               &   \theta_{2}
\end{bmatrix}\otimes\theta')\]\\
\[=(\begin{bmatrix}
  \phi_{1}\otimes \phi' &   0      \\
  0   & \phi_{2}\otimes \phi'
\end{bmatrix}, \begin{bmatrix}
  \psi_{1}\otimes \psi' &   0      \\
  0   & \psi_{2}\otimes \psi'
\end{bmatrix}, \begin{bmatrix}
  \theta_{1}\otimes \theta' &   0      \\
  0   & \theta_{2}\otimes \theta'
\end{bmatrix})\,\,\,\,\,\,\,\,\cdots (\sharp')\]
\\%
The desired identity now follows from $(\sharp)$ and $(\sharp')$.
\item The proof of this equality is similar to the foregoing proof.
\end{enumerate}
\end{proof}

\section{FURTHER PROBLEMS}
We constructed a bifunctorial operation $\overline{\otimes}_{3}$ which is such that if $X$ (respectively $Y$) is a $3-$matrix factorization of $f\in R$ (respectively $g\in S$), then $X\overline{\otimes}_{3} Y$ is a $3-$matrix factorization of $fg\in K[x_{1},x_{2},\cdots, x_{m},y_{1},y_{2},\cdots, y_{m}]$. An interesting question would be to find a counterpart to the Yoshino tensor product of matrix factorizations, that is a bifunctorial operation $\widehat{\otimes}_{3}$ which is such that if $X$ (respectively $Y$) is a $3-$matrix factorization of $f\in R$ (respectively $g\in S$), then $X\widehat{\otimes}_{3} Y$ is a $3-$matrix factorization of $f+g\in K[x_{1},x_{2},\cdots, x_{m},y_{1},y_{2},\cdots, y_{m}]$.


\bibliography{fomatati_ref}
\addcontentsline{toc}{section}
{References}
\end{document}